\documentclass[11pt,leqno]{article}

\usepackage{authblk}

\title{Implementation of the Continuous-Discontinuous Galerkin Finite Element Method}

\author[1]{Andrea Cangiani\thanks{andrea.cangiani@le.ac.uk}}
\author[2]{John Chapman\thanks{john.chapman@durham.ac.uk}}
\author[1]{Emmanuil Georgoulis\thanks{emmanuil.georgoulis@le.ac.uk}}
\author[2]{Max Jensen\thanks{m.p.j.jensen@durham.ac.uk}}

\affil[1]{Department of Mathematics, University of Leicester, University Road, Leicester, United Kingdom}
\affil[2]{Department of Mathematics, University of Durham, Durham, United Kingdom}

\date{\today}

\usepackage[paper=a4paper,dvips,top=2.0cm,left=2.0cm,right=2.0cm,
    foot=1cm,bottom=2.0cm]{geometry}

\usepackage[dvips]{color}

\usepackage{graphicx}
\usepackage{subfigure}

\usepackage[]{amsmath}

\usepackage{amsfonts}
\usepackage{amssymb}
\usepackage{textcomp}

\usepackage{amsthm}

\usepackage{amsopn}

\usepackage{xspace}

\usepackage{enumitem}

\usepackage{perpage}
\MakePerPage{footnote}

\usepackage[pagewise]{lineno} 
\newcommand*\patchAmsMathEnvironmentForLineno[1]{%
  \expandafter\let\csname old#1\expandafter\endcsname\csname #1\endcsname
  \expandafter\let\csname oldend#1\expandafter\endcsname\csname end#1\endcsname
  \renewenvironment{#1}%
     {\linenomath\csname old#1\endcsname}%
     {\csname oldend#1\endcsname\endlinenomath}}%
\newcommand*\patchBothAmsMathEnvironmentsForLineno[1]{%
  \patchAmsMathEnvironmentForLineno{#1}%
  \patchAmsMathEnvironmentForLineno{#1*}}%
\AtBeginDocument{%
\patchBothAmsMathEnvironmentsForLineno{equation}%
\patchBothAmsMathEnvironmentsForLineno{align}%
\patchBothAmsMathEnvironmentsForLineno{flalign}%
\patchBothAmsMathEnvironmentsForLineno{alignat}%
\patchBothAmsMathEnvironmentsForLineno{gather}%
\patchBothAmsMathEnvironmentsForLineno{multline}%
}

\newtheorem{theorem}[equation]{Theorem}
\newtheorem{definition}[equation]{Definition}

\numberwithin{equation}{section}
\numberwithin{figure}{section}
\numberwithin{table}{section}

\newcommand\CC{\Lang{\mbox{C++}}\xspace}

\newcommand\Lang[1]{\textsc{#1}}

\newcommand{\bdm}{\begin{displaymath}}
\newcommand{\edm}{\end{displaymath}}
\newcommand{\beq}{\begin{equation}}
\newcommand{\eeq}{\end{equation}}

\newcommand{\newVar}[2]{\newcommand{#1}{\ensuremath{#2}\xspace}}
\newVar\Naturals{\mathbb{N}}
\newVar\Integers{\mathbb{Z}}
\newVar\Rationals{\mathbb{Q}}
\newVar\Reals{\mathbb{R}}
\newVar\Complex{\mathbb{C}}

\newcommand{\code}[1]{\texttt{#1}}

\newcommand{\norm}[1]{\ensuremath{\lVert#1\rVert}}

\newcommand{\set}[1]{\ensuremath{\{#1\}}}

\newcommand{\Average}[1]{\ensuremath{\{ \! \! \{ #1 \}  \! \! \}} }
\newcommand{\average}[1]{\ensuremath{\Average{#1}}}

\usepackage{stmaryrd}
\newcommand{\Jump}[1]{\ensuremath{\llbracket #1 \rrbracket}}
\newcommand{\jump}[1]{\ensuremath{\Jump{#1}}}

\newcommand{\transpose}[0]{{\ensuremath{\top}}}
\newcommand{\ex}[0]{\ensuremath{\mathrm{e}}}
\newcommand{\Order}[1]{{\ensuremath{\mathcal{O}(#1)}}}

\newVar\pqTree{{\cal P}}
\newVar\groups{\mathbf{G}}
\newVar\lastGroup{\ell}
\newVar{\currentId}{\ensuremath{\id}\xspace}

\renewcommand\epsilon{\varepsilon}
\renewcommand\phi{\varphi}
\renewcommand\theta{\vartheta}

\newcommand{\Poly}[0]{{\ensuremath{\mathbb{P}}}}

\renewcommand{\cdot}[0]{\ensuremath{\cdotp}}
\newcommand{\half}[0]{{\ensuremath{\frac{1}{2}}}}
\usepackage{units}
\newcommand{\fhalf}[0]{{\ensuremath{\nicefrac{1}{2}}}}

\newcommand{\cdG}[0]{{\ensuremath{\textrm{cdG}}}}
\newcommand{\dG}[0]{{\ensuremath{\textrm{dG}}}}
\newcommand{\cG}[0]{{\ensuremath{\textrm{cG}}}}

\newcommand{\VdG}[1]{\ensuremath{V_\dG^{#1}}}
\newcommand{\VcdG}[1]{\ensuremath{V_\cdG^{#1}}}
\newcommand{\cdGh}[0]{{\ensuremath{\cdG,h}}}
\newcommand{\dGh}[0]{{\ensuremath{\dG,h}}}

\newcommand{\cont}[0]{\ensuremath{\cG}}
\newcommand{\discont}[0]{\ensuremath{\dG}}

\newcommand{\Th}[0]{{\ensuremath{\mathcal{T}_h}}}

\newcommand{\TD}[0]{{\ensuremath{\mathcal{T}_h^{\discont}}}}
\newcommand{\TC}[0]{{\ensuremath{\mathcal{T}_h^{\cont}}}}

\newcommand{\EinTh}[0]{{\ensuremath{E \in \Th}}}

\newcommand{\Eho}[0]{\ensuremath{{\mathcal{E}_h^o}}}

\newcommand{\Eh}[0]{\ensuremath{{\mathcal{E}_h}}}
\newcommand{\EhC}[0]{\ensuremath{{\mathcal{E}_h^\cont}}}
\newcommand{\EhD}[0]{\ensuremath{{\mathcal{E}_h^\discont}}}
\newcommand{\GC}[0]{\ensuremath{{\Gamma^\cont}}}
\newcommand{\GD}[0]{\ensuremath{{\Gamma^\discont}}}
\newcommand{\Gout}[0]{\ensuremath{\Gamma^{\text{out}}}}
\newcommand{\Gin}[0]{\ensuremath{\Gamma^{\text{in}}}}
\newcommand{\einEh}[0]{{\ensuremath{e \in \Eh}}}
\newcommand{\einEho}[0]{{\ensuremath{e \in \Eho}}}
\newcommand{\einEhC}[0]{{\ensuremath{e \in \EhC}}}

\newcommand{\pO}[0]{{\ensuremath{\partial \Omega}}}
\newcommand{\OD}[0]{\ensuremath{\Omega^\discont}}
\newcommand{\OC}[0]{\ensuremath{\Omega^\cont}}
\newcommand{\pOC}[0]{\ensuremath{\partial \OC}}

\renewcommand{\div}[0]{{\ensuremath{\mathrm{div}}}}

\newcommand{\SH}[3]{{\ensuremath{H^{#1}_{#2}(#3)}}}

\newcommand{\SL}[2]{{\ensuremath{L^{#1}(#2)}}}

\newcommand{\fenothing}[0]{{\code{FE\_NOTHING}~}}
\newcommand{\dealii}[0]{{\code{deal.ii}~}}

\begin{document}

\maketitle

\begin{abstract}
For the stationary advection-diffusion problem the standard continuous Galerkin method is unstable without some additional control on the mesh or method. The interior penalty discontinuous Galerkin method is stable but at the expense of an increased number of degrees of freedom. The hybrid method proposed in \cite{CGJ06} combines the computational complexity of the continuous method with the stability of the discontinuous method without a significant increase in degrees of freedom. We discuss the implementation of this method using the finite element library \dealii and present some numerical experiments.
\end{abstract}

\section{Introduction} \label{sec:Introduction}

We consider the advection-diffusion equation
\begin{align} 
		-\epsilon \Delta u + b \cdot \nabla u & = f  \qquad \mathrm{in}~ \Omega \subset \Reals^d \label{ADeqn1}\\
		u & = g   \qquad \mathrm{on}~ \pO \label{ADeqn2}
\end{align}
with $0 < \epsilon \ll 1$, $b \in W^\infty(\div, \Omega)$, $f \in \SL{2}{\Omega}$ and $g \in \SH{\fhalf}{}{\Omega}$. For simplicity we assume the region $\Omega$ is polygonal. We also assume $\rho := -\half \nabla \cdot b \ge 0$ and then we have a weak solution $u \in \SH{1}{}{\Omega}$.

It is well known that this problem can exhibit boundary or internal layers in the convection dominated regime and that for the standard continuous Galerkin (cG) formulation these layers cause non-physical oscillations in the numerical solution. Several adaptations to the cG method are effective but space does not allow their discussion here. We refer readers to \cite{RST08} for a full description of these approaches. Discontinuous Galerkin (dG) methods also offer a stable approach for approximating this problem. However the number of degrees of freedom required for dG methods is in general considerably larger than for cG methods.

We describe an alternative approach also studied in \cite{CGJ06,DP02,DFG07}. A dG method is applied on the layers and a cG method away from the layers. We call this approach the continuous-discontinuous Galerkin (cdG) method. The hypothesis is that provided the layers are entirely contained in the dG region the instability they cause will not propagate to the cG region. Note that in our formulation there are no transmission conditions at the join of the two regions.

Here we present the cdG method and discuss its implementation using the \dealii finite element library. We additionally provide some numerical experiments to highlight the performance of the method.

\section{Finite Element Formulation} \label{sec:FEFormulation}

Assume that we can identify a decomposition of $\Omega := \OC \cup \OD$ where it is appropriate to apply the cG and dG methods respectively. We do not consider specific procedures to achieve this here, but generally it will be that we wish all boundary and internal layers to be within $\OD$. Identifying these regions can be done a priori in some cases or a posteriori based on the solution of a dG finite element method. Consider a triangulation $\Th$ of $\Omega$ which is split into two regions $\TC$ and $\TD$ where we will apply the cG and dG methods respectively. For simplicity we assume that the regions $\TC$ and $\TD$ are aligned with the regions $\OC$ and $\OD$ and the set $J$ contains edges which lie in the intersection of the two regions. Call the mesh skeleton $\Eh$ and the internal skeleton $\Eho$. Define $\Gamma$ as the union of boundary edges and the inflow and outflow boundaries by
\begin{align*}
		\Gin = & \set{x \in \pO : b \cdot n \le 0} \\
		\Gout = & \set{x \in \pO : b \cdot n > 0} 
\end{align*}
where $n$ is the outward pointing normal. Define $\GC$ (resp. $\GD$) to be the intersection of $\Gamma$ with $\TC$ (resp. $\TD$). By convention we say that the edges of $J$ are part of the discontinuous skeleton $\EhD$ and $\EhC := \Eh \setminus \EhD$. With this convention there is potentially a discontinuity of the numerical solution at $J$. Elements of the mesh are denoted $E$, edges (resp. faces in 3d) by $e$ and denote by $h_E$ and $h_e$ the diameter of an element and an edge, defined in the usual way. 

The jump $\jump{\,\cdot}$ and average $\average{\,\cdot}$ of a scalar or vector function on the edges in $\Eh$ are defined as in, e.g., \cite{ABCM01}.


\begin{definition}
  Define the cdG space to be
  \begin{equation} \label{cdGSpace}
  \begin{split}  
    \VcdG{k} &:= \set{ v \in \SL{2}{\Omega} : v|_E \in \Poly^k, v|_{\pO \cap \pOC} = g, v|_{\OC} \in \SH{1}{}{\OC}}
  \end{split}
  \end{equation}
  where $\Poly^{k}$ is the space of polynomials of degree at most $k$ supported on $E$. This is equivalent to applying the usual cG space on $\OC$ and a dG space on $\OD$.
\end{definition}
  
We may now define the interior penalty cdG method: Find $u_h \in \VcdG{k}$ such that for all $v_h \in \VcdG{k}$
\begin{equation}
	B_\epsilon (u_h,v_h) = B_d(u_h,v_h) + B_a(u_h,v_h) = L_\epsilon(f,g;v_h) \nonumber
\end{equation}
where
\begin{equation}
	\begin{split}
		B_d(u_h,v_h) &=  \sum_\EinTh \left[ \int_E \epsilon \nabla u_h \cdot \nabla v_h - \int_E (b \cdot \nabla v_h) u_h - \int_E (\nabla\cdot b)u_h v_h \right] \nonumber \\
			& \qquad +\sum_\einEh \left[ \int_e \sigma \frac{\epsilon}{h_e} \jump{u_h} \cdot \jump{v_h} - \int_e \left( \average{\epsilon \nabla u_h}\cdot \jump{v_h} + \theta \average{\epsilon \nabla v_h}\cdot \jump{u_h}\right) \right] \nonumber \\
		B_a(u_h,v_h) &= \sum_\einEho \int_e b \cdot \jump{v_h} u_h^- + \sum_{e \in \Gout} \int_e (b \cdot n) u_h v_h \nonumber
	\end{split}
\end{equation}
and
\begin{equation}
	L_\epsilon (f,g;v_h) = \sum_\EinTh \int_E fv_h + \sum_{e \in \Gamma} \left[ \int_e \left(\sigma \frac{\epsilon}{h_e} v_h - \theta \epsilon \nabla \cdot v_h \right) g \right] - \sum_{e \in \Gin} \int_e (b \cdot n)v_hg. \nonumber
\end{equation}
Here $\sigma$ is the penalization parameter and $\theta \in \set{-1,0,1}$. Note that through the definition of $\VcdG{k}$ the edge terms are zero on $\EhC$ and the method reduces to the standard cG FEM. If we take $\Th = \TD$, i.e., the entire triangulation as discontinuous, we get the interior penalty (IP) family of dG FEMs (see, e.g., \cite{ABCM01}).


The work of \cite{LN00} shows that the cG method is the limit of the dG method as $\sigma \rightarrow \infty$. A reasonable hypothesis is that the solution to the cdG method is the limit of the solutions to the dG method as the penalty parameter $\sigma \rightarrow \infty$ on $\einEhC$, i.e., super penalising the edges in $\EhC$. Call $\sigma_\cont$ and $\sigma_\discont$ the penalty parameters for edges in $\EhC$ and $\EhD$ respectively. Call the numerical solution for the cdG problem $u_\cdGh \in \VcdG{k}$. The solution to the pure dG problem on the same mesh is denoted $u_\dGh \in \VdG{k}$ where $\VdG{k}$ is the usual piecewise discontinuous polynomial space on $\Th$. Then we have
\begin{theorem} \label{thm:convergence}
  The dG solution converges to the cdG solution of \eqref{ADeqn1}-\eqref{ADeqn2} as $\sigma_\cont \rightarrow \infty$, i.e.,
  \begin{equation*}
    \lim_{\sigma_\cont \rightarrow \infty} (u_\cdGh - u_\dGh) = 0
  \end{equation*}
\end{theorem}
We do not prove this result here but direct readers to \cite{CCGJ12} for a full discussion. Although this result does not imply stability of the cdG method (indeed, for the case where the $\OC$ region is taken to be the whole of $\Omega$ it shows that the cdG method has the same problems as the cG method), but does indicate that investigation of the cdG method as an intermediate stage between cG and dG is justified. Hence, it aids into building an understanding of the convergence and stability properties of the cdG method, based on what is known for cG and dG. This, in turn, is of interest, as the cdG method offers substantial reduction in the degrees of freedom of the method compared to dG.

\section{Numerical Implementation} \label{sec:Implementation}

The cdG method poses several difficulties in implementation. One approach is to use the super penalty result of Theorem \ref{thm:convergence} to get a good approximation to the cdG solution. However this will give a method with the same number of degrees of freedom as dG. We therefore present an approach to implement the cdG method with the appropriate finite element structure. We discuss this approach with particular reference to the \dealii finite element library \cite{BHK07,DealIIRef}. This is an open source \CC library designed to streamline the creation of finite element codes and give straightforward access to algorithms and data structures. We also present some numerical experiments.

\subsection{Implementation in \dealii} \label{sec:Implementationdealii}

The main difficulty in implementing a cdG method in \dealii is the understandable lack of a native cdG element type. In order to assign degrees of freedom to a mesh in \dealii the code must be initialised with a \code{Triangulation} and then instructed to use a particular finite element basis to place the degrees of freedom. Although it is possible to initialise a \code{Triangulation} with the dG and cG regions set via the \code{material\_id} flag, no appropriate element exists. In the existing \dealii framework it would be difficult to code an element with the appropriate properties. A far more robust approach is to use the existing capabilities of the library and therefore allow access to other features of \dealii\!. For instance without the correct distribution of degrees of freedom the resulting sparsity pattern of the finite element matrix would be suboptimal, i.e., containing more entries than required by the theory and therefore reducing the benefit of shrinking the number of degrees of freedom relative to a dG method.

The \dealii library has the capability to handle problems with multiple equations applied to a single mesh such as the case of a elastic solid fluid interaction problem. In our case we wish to apply different methods to the same equation on different regions of the mesh, which is conceptually the same problem in the \dealii framework. In addition we will use the $hp$ capability of the library.

The \dealii library has the capability to create collections of finite elements, \code{hp::FECollection}. Here multiple finite elements are grouped into one data structure. As the syntax suggests the usual use is for $hp$ refinement to create a set of finite elements of the same type (e.g., scalar Lagrange elements \code{FE\_Q} or discontinuous elements \code{FE\_DGQ}) of varying degree. Unfortunately it is not sufficient to create a \code{hp::FECollection} of cG and dG elements as the interface between the two regions will still be undefined. In order to create an admissible collection of finite elements we use \fenothing\!. This is a finite element type in \dealii with zero degrees of freedom. Using the $\code{FESystem}$ class we create two vector-valued finite element types $(\code{FE\_Q},\fenothing\!)$ and $(\fenothing\!,\code{FE\_DGQ})$ and combine them in a \code{hp::FECollection}. We apply the first $\code{FESystem}$ on the cG region, and the second on the dG region. Now when we create a \code{Triangulation} initialised with the location of cG and dG elements the degrees of freedom can be correctly distributed according to the finite element defined by \code{hp::FECollection}.

When assembling the matrix for the finite element method we need only be careful that we are using the correct element of \code{hp::FECollection} and the correct part of $\code{FESystem}$. The most difficult case is on the boundary $J$ where from a dG element we must evaluate the contribution from the neighbouring cG element (note that in the cdG method a jump is permissible on $J$).

If we implement the cdG method in \dealii in this way we create two solutions: one for the \code{FE\_Q}-\fenothing component and another for the \fenothing-\code{FE\_DGQ} component. Consider a domain $\Omega = (0,1)^2$ in $\Reals^2$, $b = (1,1)^\transpose$ and $\theta = -1$. The Dirichlet boundary conditions and the forcing function $f$  are chosen so that the analytical solution is
\begin{equation} \label{Ex1TrueSol}
	u(x,y)=x+y(1-x)+\frac{\ex^{-\frac{1}{\epsilon}}-\ex^{-\frac{(1-x)(1-y)}{\epsilon}}}{1-\ex^{-\frac{1}{\epsilon}}}.
\end{equation}
This solution exhibits an exponential boundary layer along $x=1$ and $y=1$ of width $\Order{\epsilon}$, $\epsilon = 10^{-6}$. We solve the finite element problem on a 1024 element grid and fix $\OC = [0, 0.707)^2$. This is larger than is required for stability (see Example 1 below) but shows the behaviour of \fenothing more clearly. We show each of the components of \code{FE\_SYSTEM} and the combined solution. For comparison we also show the dG finite element solution for the same problem.

\begin{figure}[phtb]
	\centering
	\subfigure[cG-\fenothing]{\includegraphics[width=0.30\textwidth]{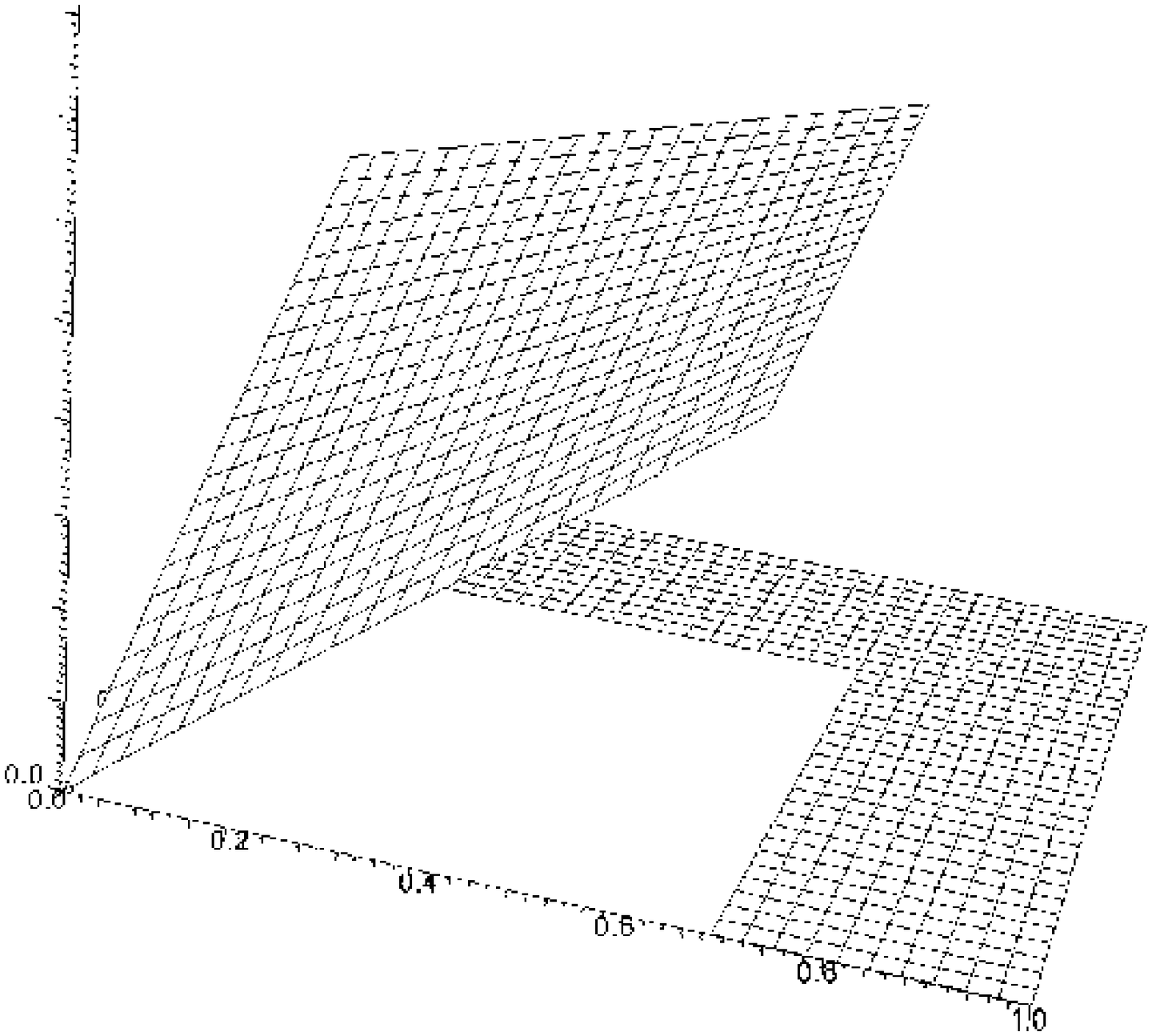}} \hspace{20pt}
	\subfigure[\fenothing\!-dG]{\includegraphics[width=0.30\textwidth]{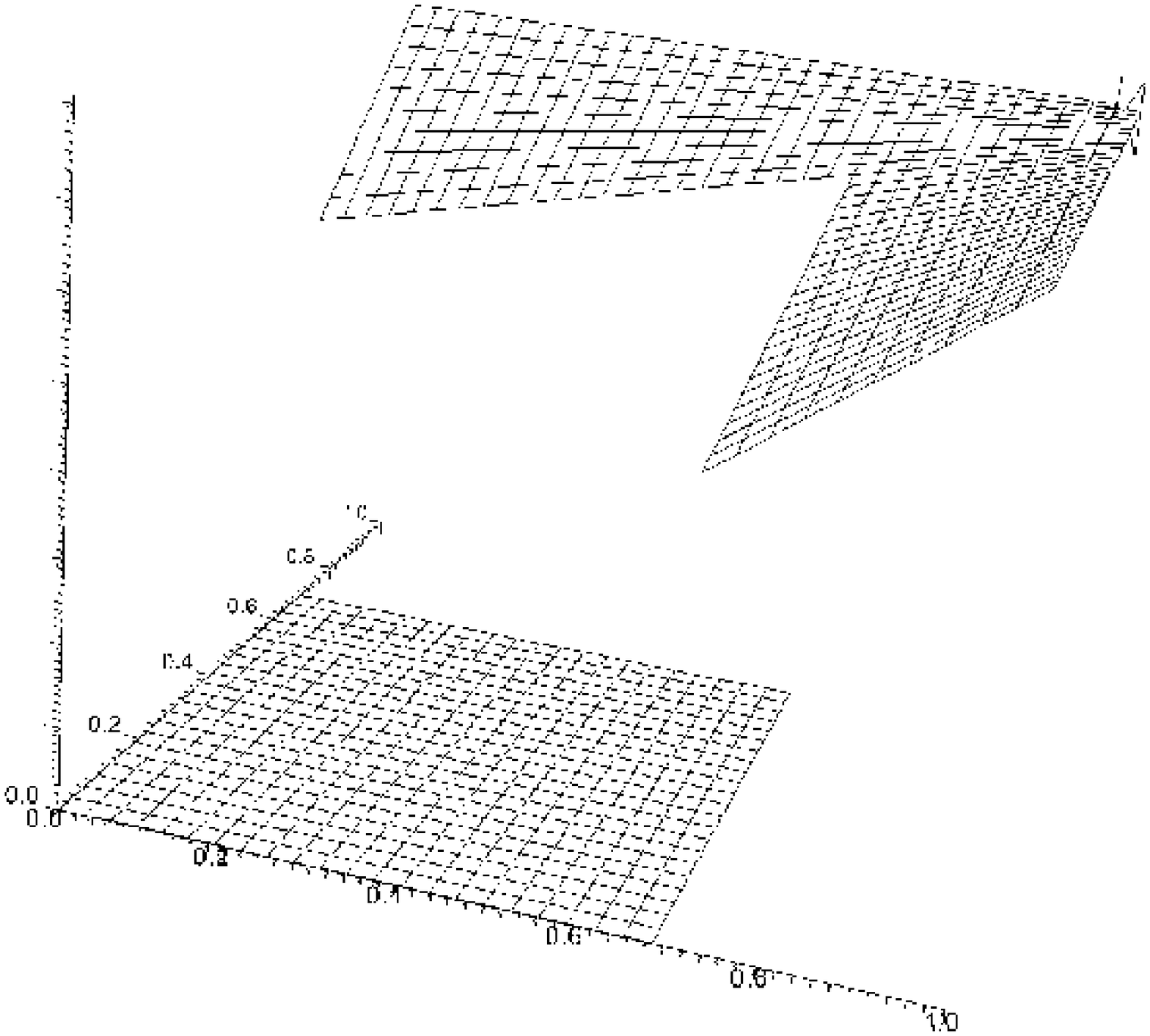}}  \\
	\subfigure[Combined cdG solution]{\includegraphics[width=0.30\textwidth]{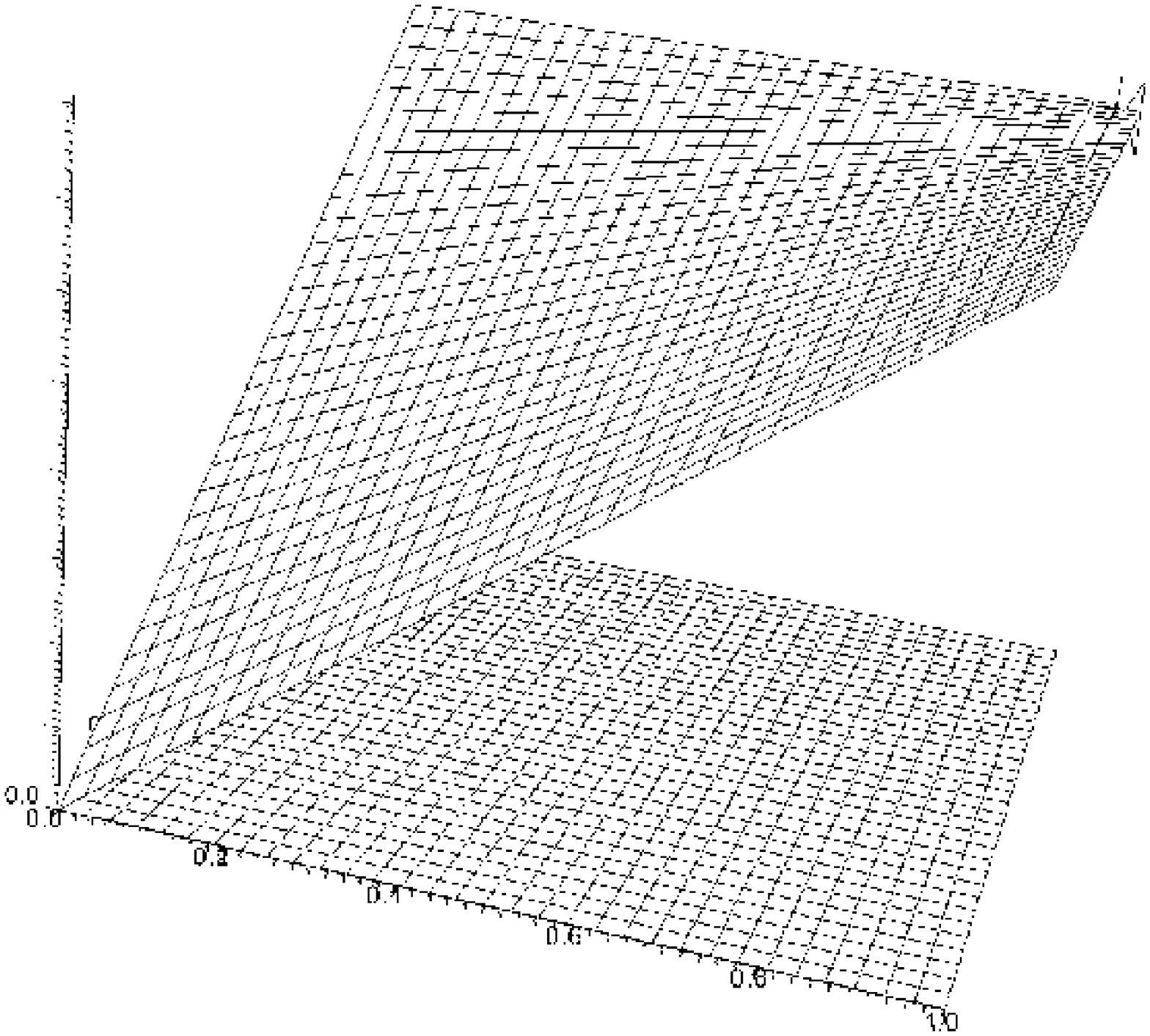}} \hspace{20pt}	
	\subfigure[dG solution]{\includegraphics[width=0.30\textwidth]{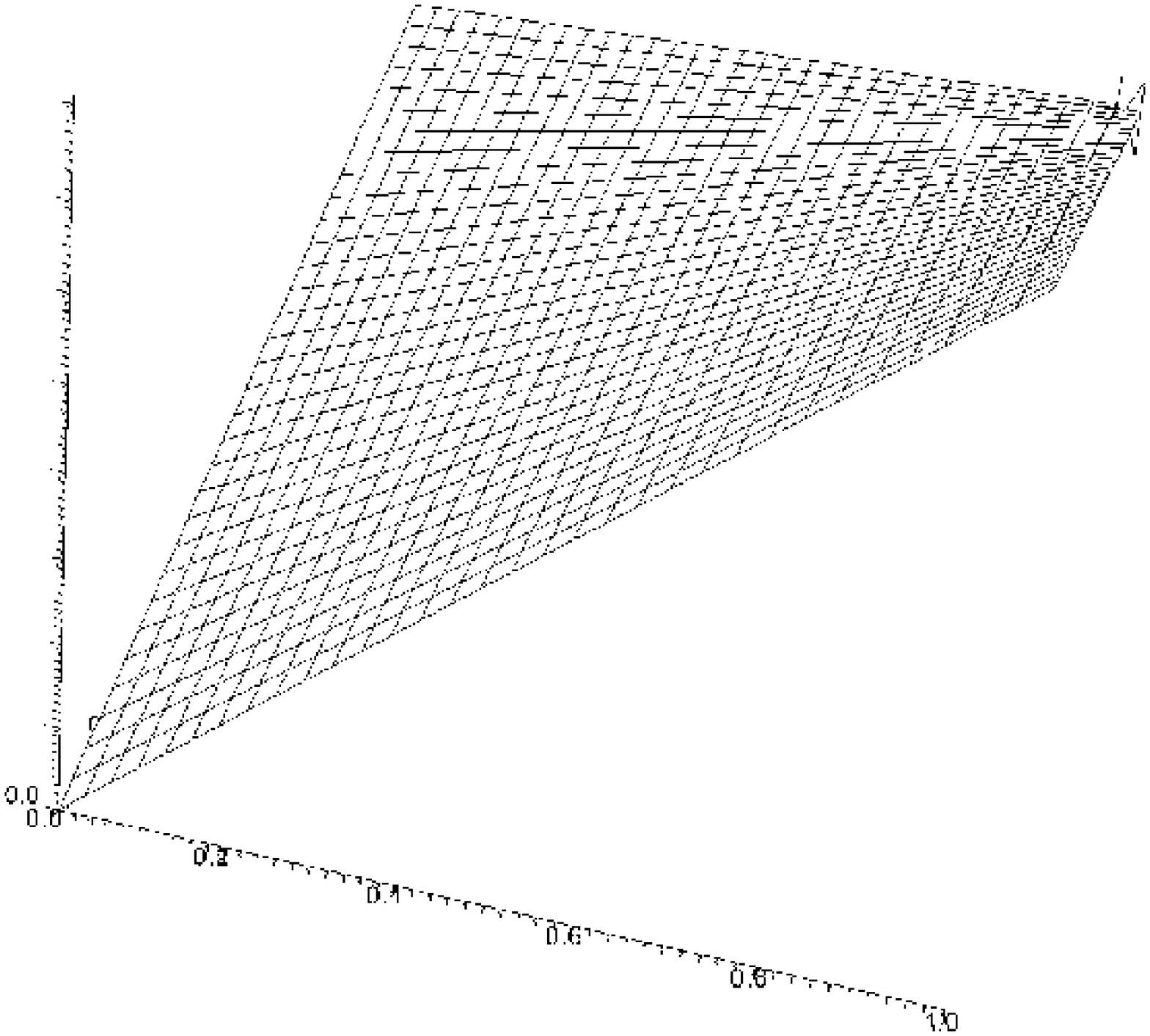}} 	
	\caption{Solution components of \fenothing implementation applied to Example 1 with $\epsilon = 10^{-6}$.}
	\label{fig:FENothing}
\end{figure}

One advantage of following the \dealii framework is that the data structures will allow the implementation of $hp$ methods. In fact we can envisage the implementation of a $hpe$ method where at each refinement there is the possibility to change the mesh size, polynomial degree or the element type. We propose no specific scheme here but simply remark that implementing a $hpe$ method is relatively straightforward with the \fenothing approach. 

\subsection{Numerical Examples} \label{sec:NumericalExamples}

We present two numerical experiments highlighting the performance of the cdG method. Both examples present layers when $\epsilon$ is small enough. In each case we fix the region where the continuous method is to be applied then vary $\epsilon$. This causes the layer to steepen. In the advection dominated regime, i.e., $\epsilon$ large and no steep layer present, we see the cdG solution approximates the true solution well. As we make $\epsilon$ smaller the layer forms and extends into the continuous region. As $\epsilon$ becomes smaller still the layer leaves the continuous region and the performance of the dG and cdG method is indistinguishable. In each experiment we pick the $\OC$ and $\OD$ regions so that with the given refinement the region $\TD$ consists of exactly one layer of elements and coincides with $\OD$. 

\paragraph{Example 1}

Consider again the problem with true solution \eqref{Ex1TrueSol} presented above. We solve the finite element problem on a 1024 element grid and fix $\OC = [0, 0.96875)^2$ so exactly one row of elements is in $\OD$. As we vary $\epsilon$ the layer sharpens and moves entirely into the dG region.

\begin{figure}[phtb]
	\centering
	\subfigure[$\norm{u_\cdGh -  u_\dGh}_\SL{2}{\Omega}$]{\includegraphics[width=0.45\textwidth]{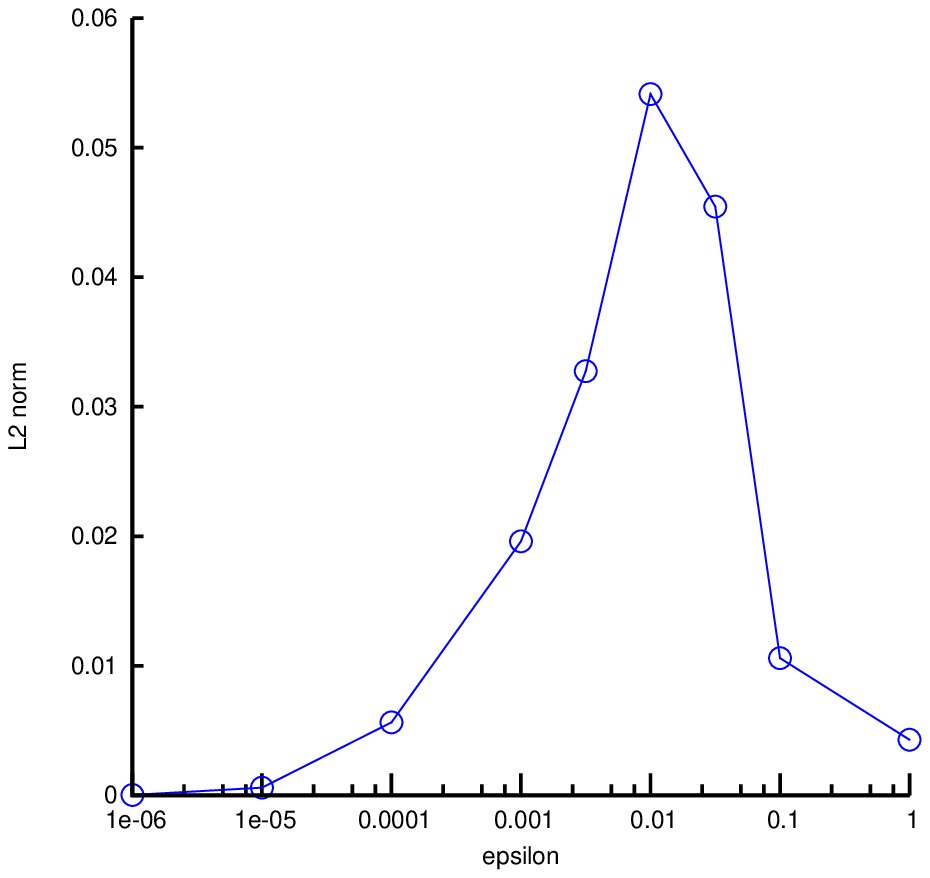}} \hspace{20pt}
	\subfigure[$\norm{u_\cdGh -  u_\dGh}_\SL{\infty}{\Omega}$]{\includegraphics[width=0.45\textwidth]{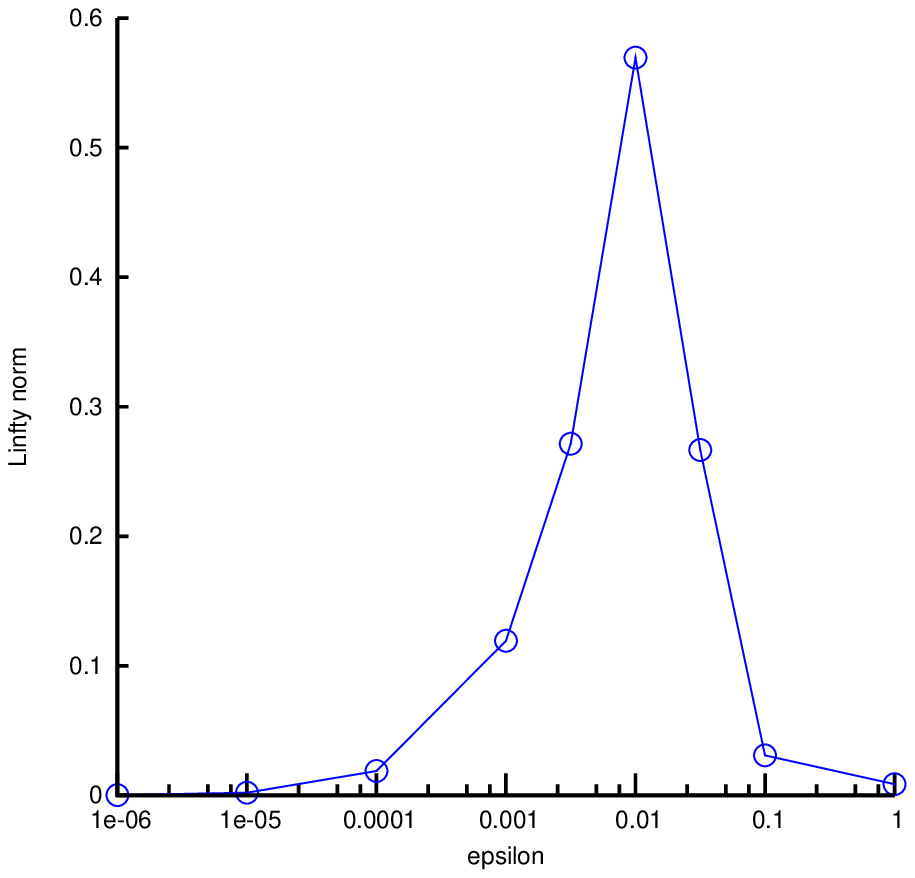}}\\
	\caption{Decreasing $\epsilon$ with a fixed $\Omega$ decomposition in Example 1.The maximum difference in either norm occurs when the layer is sharp but not contained entirely in $\OD$.}
	\label{fig:Example1A}
\end{figure}

As we can see from Figure \ref{fig:Example1A} before the layer has formed the two methods perform well. As the layer begins to form with decreasing $\epsilon$ it is not entirely contained in the discontinuous region and the error peaks. As the layer sharpens further it is entirely contained in the discontinuous region and the difference between the two solutions becomes negligible.

\paragraph{Example 2}

Now we look at a problem with an internal layer. Let the advection coefficient be given by $b=(-x,y)^\transpose$ and pick the boundary conditions and right hand side $f$ so that the true solution is
\begin{equation}
	u(x,y) = (1-y^2)\textrm{erf}\left(\frac{x}{\sqrt{2\epsilon}}\right), \nonumber
\end{equation}
where $\textrm{erf}$ is the error function defined by
\begin{equation}
	\textrm{erf}(x) = \frac{2}{\sqrt{\pi}} \int_0^x \ex^{-t^2} \mathrm{d}t. \nonumber
\end{equation}

We solve on the region $\Omega = (-1,1)^2$. The solution has an internal layer along $y=0$ of width $\Order{\sqrt{\epsilon}}$ and we fix $\OC = \set{(x,y) : x \in [-1,-0.0625)\cup(0.0625,1], y \in [-1,1]}$.

\begin{figure}[phtb]
	\centering
	\subfigure[$\epsilon = 10^{-2}$]{\includegraphics[width=0.32\textwidth]{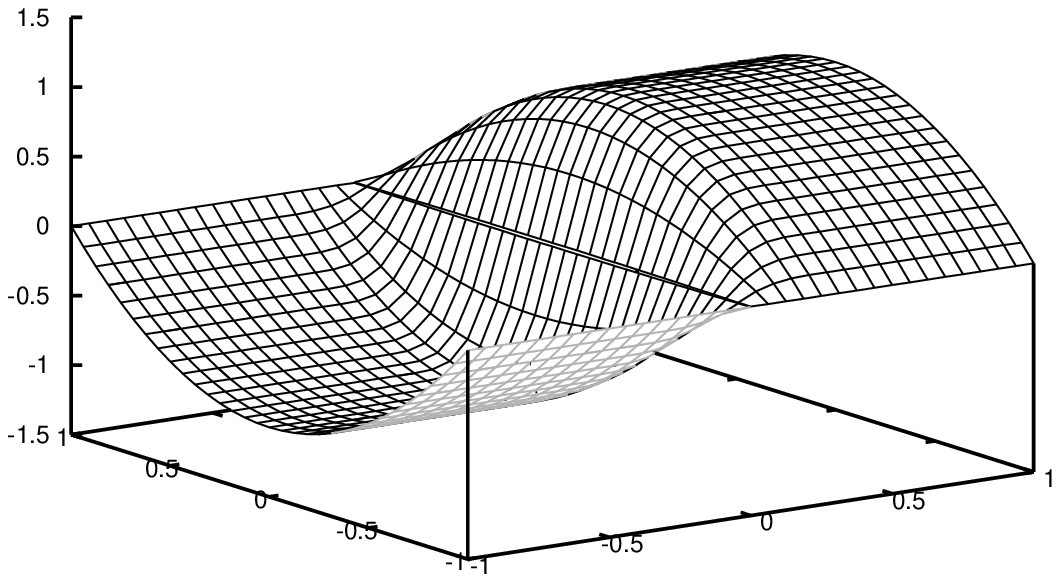}} \hspace{1pt}
	\subfigure[$\epsilon = 10^{-4}$]{\includegraphics[width=0.32\textwidth]{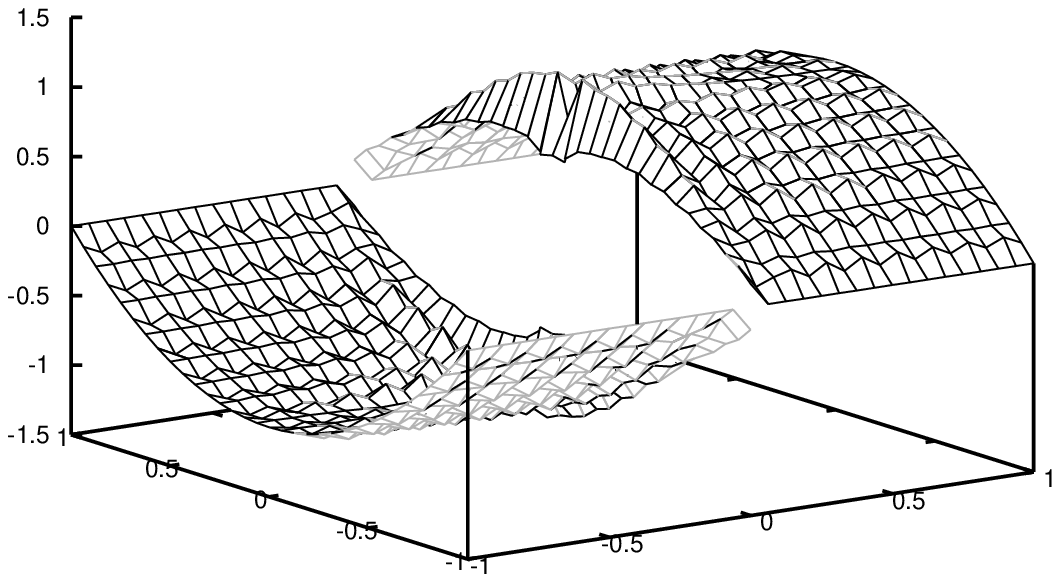}} \hspace{1pt}
	\subfigure[$\epsilon = 10^{-6}$]{\includegraphics[width=0.32\textwidth]{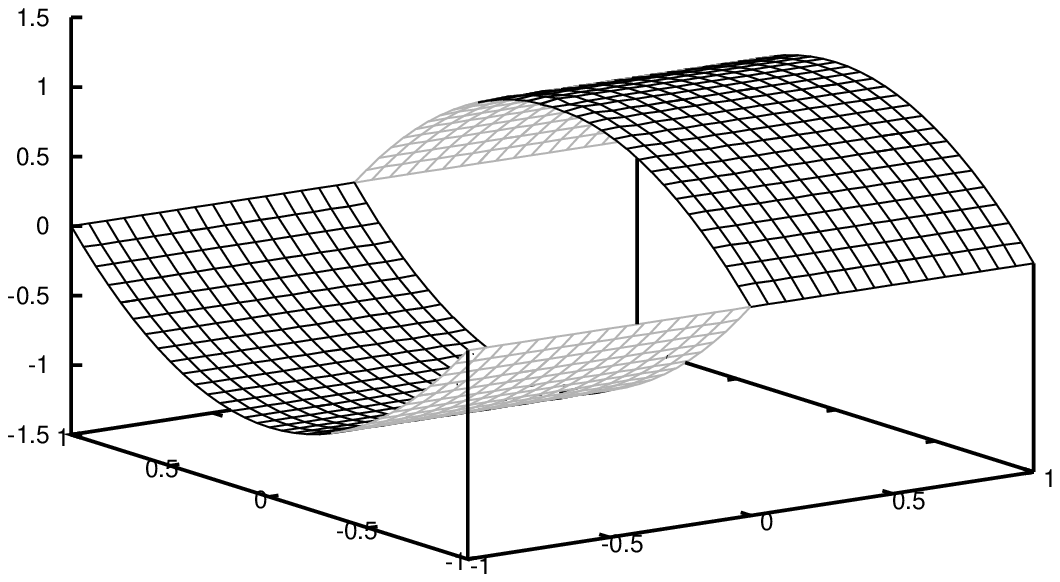}}
	\caption{The cdG solutions for Example 2 for various $\epsilon$. When $\epsilon = 10^{-4}$ the layer is steep enough to cause oscillations but not sharp enough to be contained entirely in $\OD$. In this case the oscillations are clearly visible, but they are not present when $\epsilon = 10^{-6}$ as the layer has moved entirely within $\OD$.}
	\label{fig:Example2B}
\end{figure}

\begin{figure}[phtb]
	\centering
	\subfigure[$\norm{u_\cdGh -  u_\dGh}_\SL{2}{\Omega}$]{\includegraphics[width=0.45\textwidth]{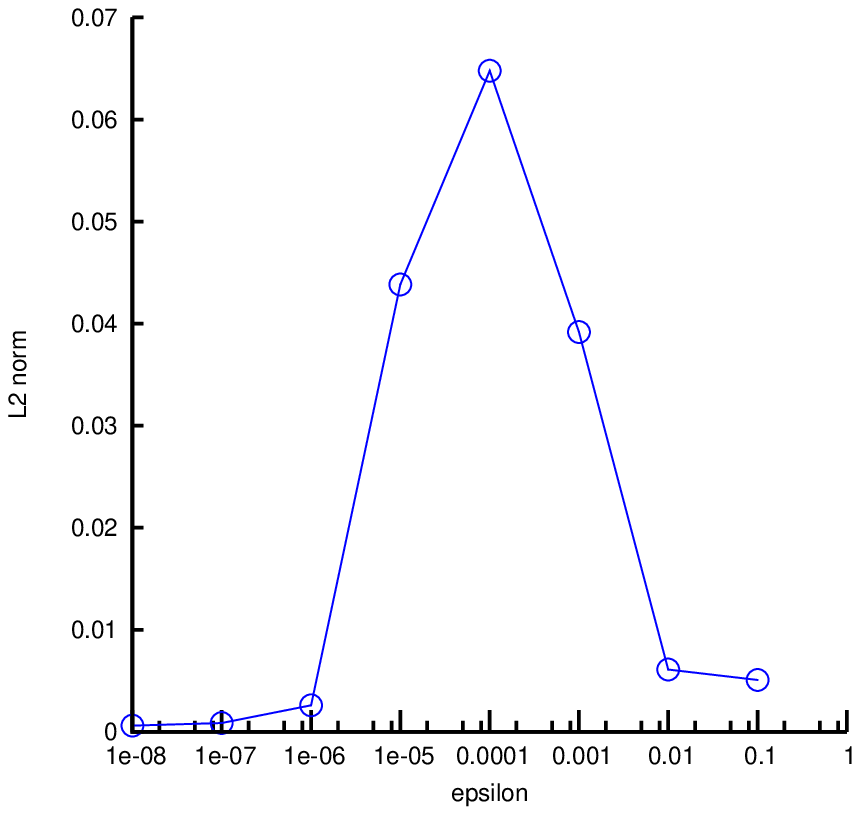}} \hspace{20pt}
	\subfigure[$\norm{u_\cdGh -  u_\dGh}_\SL{\infty}{\Omega}$]{\includegraphics[width=0.45\textwidth]{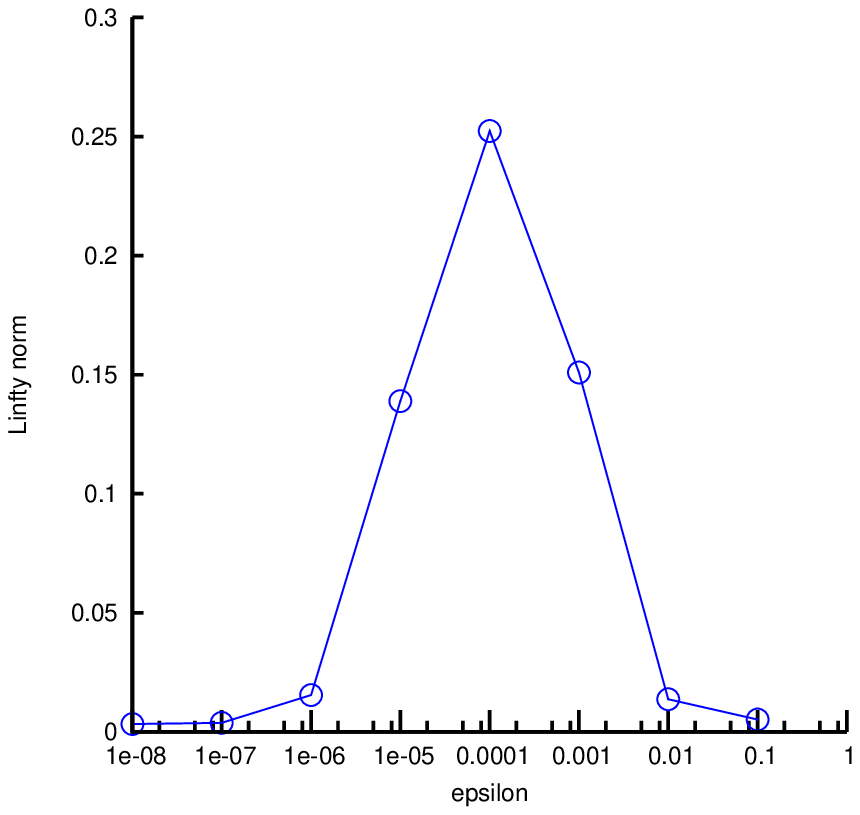}}\\
	\caption{Decreasing $\epsilon$ with a fixed $\Omega$ decomposition in Example 2. As in Figure \ref{fig:Example1A} the maximum difference in either norm occurs when the layer is sharp but not contained entirely in $\OD$.}
	\label{fig:Example2A}
\end{figure}

In Figure \ref{fig:Example2A} we notice same the same behaviour as in Example 1. If the layer exists it must be contained within the discontinuous region for the two methods to perform equivalently. In Figure \ref{fig:Example2B} we can see the cdG solution for various $\epsilon$ with the oscillations clearly visible when $\epsilon = 10^{-4}$. When the layer is sharpened, the oscillations disappear.

\end{document}